# ON DIVISIBILITY PROPERTIES OF TRUNCATED BINOMIAL SERIES

Anatoly A. Grinberg

anatoly_gr@yahoo.com

## ABSTRACT

The divisibility of truncated binomial series by their exponent n is analyzed. Divisibility is shown to depends on the divisibility characteristics of the integers constituting the binomials. Series division by the highest possible powers of n is examined.

In number theory problems of compatibility, evenness, coprime status and divisibility of equation terms arise often. In particular, divisibility of binomial series by their exponent n is an interesting problem[1]. Here we consider the divisibility of truncated binomial series. The divisibility of such series strongly depends on the divisibility characteristics of the binomial's constituent integers. Our goal is to find the maximum exponent power of n for the divisibility.

It will be assumed that the greatest common divisor (GCD) of these numbers has been extracted. Everywhere we operate only with integer numbers.

## I. Truncated Newton's binomial series.

We consider the truncated binomial series (TBS) that is obtained by removing the side terms ($a^n$ and $b^n$) from Newton's binomial $(a+b)^n = q^n$ series. Using the notations $q = a + b$ and $Q = a^n + b^n$, the TBS can be represented by the series

$$U(a,b) = q^n - Q = \sum_{v=1}^{n-1} \binom{n}{v} a^v b^{n-v} \qquad (I.1)$$

where the symbol $\binom{n}{v}$ denotes the binomial coefficients.

In order to describe a variety of cases we represent numbers $a$ and $b$ in the form

$$a = g_a n + r_a, \quad b = g_b n + r_b, \qquad (I.2)$$

where $r_a$ and $r_b$ are the respective remainders of the division of $a$ and $b$ by n. We will assume that $g_a$, $g_b$ are not divisible by n. The case where both $r_a$ and $r_b$ equal zero must be excluded since the GCD of $a$ and $b$ was already removed.

Substituting (I.2) into (I.1) we obtain the following formula for $U(a,b)$

$$U(a,b) = (Gn + R)^n - (g_a n + r_a)^n - (g_b n + r_b)^n =$$

$$n^n \sum_{v=1}^{n-1} \binom{n}{v} g_a^v g_b^{n-v} + \sum_{v=1}^{n-1} \binom{n}{v} ([G^v R^{n-v} - g_a^v r_a^{n-v} - g_b^v r_b^{n-v}] n^v +$$

$$\sum_{v=1}^{n-1} \binom{n}{v} r_a^v r_b^{n-v} \qquad (I.3)$$

where $G = g_a + g_b$, $R = r_a + r_b$

At this point, it is convenient to introduce the following symbols

$$[a, n^k] = 1, \quad \text{if } a \text{ is divisible by } n^k, \qquad (I.4a)$$

$$[a, n^k] = 0 \qquad \text{otherwise} \qquad (I.4b)$$

The three cases will be analyzed below:

1. $[a, n] = 0$, $[b, n] = 1$,

In this case $0 < r_a < n$, $r_b = 0$ and equation (I.3) acquires the form

$$U(a,b) = n^n \sum_{v=1}^{n-1} \binom{n}{v} g_a^v g_b^{n-v} + \sum_{v=3}^{n-1} \binom{n}{v} [(g_a + g_b)^v - g_a^v] n^v$$

$$+ \left(\frac{g_b}{2}\right)(2g_a + g_b) n^3 (n-1) + r_a^{n-1} g_b n^2 \qquad (I.5)$$

Thus $U(a,b)$ is divisible by $n^2$. This can be summarized in the form

$$[U(a,b), n^2] = 1 \quad \text{if} \quad [a, n] = 1 \quad \text{or} \quad [b, n] = 1 \qquad (I.6)$$

2. $[a, n] = 0$, $[b, n] = 0$, $[(a+b), n] = 1$

In this case $0 < r_a, r_b < n$, $R = r_a + r_b = n$.

The last sum in Eq.(I.3) can be rewritten using substitution $r_b = n - r_a$

$$\sum_{v=1}^{n-1}\binom{n}{v}r_a^v r_b^{n-v} = (r_a + r_b)^n - r_a^n - (n - r_a)^n =$$

$$- r_a^n(1 + (-1)^n) - \sum_{v=1}^{n-1}\binom{n}{v}n_a^v(-r_a)^{n-v} \qquad (I.7)$$

For even $n$ Eq.(I.7) shows that $U(a,b)$ is divisible by $n$, only if $[a^n, n] = 1$. For odd $n$, the first term of the sum in (I.7) is equal to $n^2(-r_a)^{n-1}$ and $U(a,b)$ is divisible by $n^2$. Thus, we have

$$[U(a,b), n] = 1 \quad \text{for even } n \text{ if } [[a^n, n] = 1 \qquad (I.8a)$$

$$[U(a,b), n^2] = 1 \quad \text{for odd } n \qquad (I.8a)$$

3. $[a, n] = 0$, $[b, n] = 0$, $[(a+b), n] = 0$

In this case $0 < r_a, r_b < n \quad r_a + r_b \ne n$

In this general case, the integers and their algebraic combination do not contain $n$-multipliers. For prime $n$, the last sum in (I.3) is divisible by the first power of $n$ due to binomial coefficients $\binom{n}{v}$. This can be expressed using little Fermat's theorem in the form

$$\sum_{v=1}^{n-1}\binom{n}{v}r_a^v r_b^{n-v} = (r_a + r_b)^n - r_a^n - r_b^n = (\mu_{a+b} - \mu_a - \mu_b)n \qquad (I.9)$$

Divisibility of numbers $\mu_{a+b}, \mu_a, \mu_b$ (known as Fermat's quotations [2]), by $n$ is the subject of intense study in the literature [3-5]. The probability of divisibility of each of this number in the indicated above intervals are small, and probability of their coincidence to provide divisibility of shown in (I.9) algebraic combination of $\mu$ is much more smaller [1]. Thus it seems to us that in fact $U(a,b)$ can be considered as undivisible by $n^2$.

---

[1] Author express gratitude to I.D. Shkredov for helpful advise on these subject.

For non-prime n we also do not have any way of predicting the divisibility of U(a,b) by the first degree n. Thus it can be concluded that

$$[U(a,b), n] = 1 \qquad \text{for prime } n \qquad (I.10a)$$

$$[U(a,b), n] = 0 \qquad \text{for non-prime } n \qquad (I.10b)$$

## II. Divisibility of truncated trinominal series.

When analysing trinominals composed of three integers

$$q = a + b + c \qquad (II.1)$$

it is convenient to first expand $q^n$ as the sum of two numbers, (a+b) and c, then to expand the binomial $(a+b)^n$. Using this procedure and removing terms of the n-th power of a, b and c, we obtain the truncated trinomial

$$U(a, b, c) = \sum_{v=1}^{n-1} \binom{n}{v} a^v b^{n-v} + \sum_{v=1}^{n-1} \binom{n}{v} (a+b)^v c^{n-v} \qquad (II.2)$$

as the sum of two truncated binomials $U(a, b)$ and $U((a+b), c)$. The divisibility of these binomials was described in section I. The results of this analysis will be applied here. The number c will be presented in the form $c = g_c n + r_c$, similar to (I.2). Again, as in section I, we will assume that the GCD is removed from a, b and c, thereby, excluding the case $r_a = r_b = r_c = 0$.

Now we will consider the following two cases.

1. $[a, n] = 0 \quad [b, n] = 0, \quad (a+b), n] = 0, \quad [c, n] = 1$

In this case $0 < r_a, r_b < n, \quad (r_a + r_b) \neq n, \quad r_c = 0$

Under these conditions, the U(a,b) is proportional to the first degree of n, if n is a prime number (see (I.10a)), while $U(a+b, c)$ is a multiple of $n^2$ (see (I.6)). Thus, the sum of these truncated binomials leads to the linear dependence of $U(a, b, c)$ on n.

$$[U(a, b, c), n] = 1 \quad \text{if } n \text{ is prime} \quad (II.3)$$

2. $[a, n] = 0 \quad [b, n] = 0, \quad (a + b), n] = 1, \quad [c, n] = 1$

In this case we have

$$0 < r_a, r_b < n, \quad (r_a + r_b) = n, \quad r_c = 0$$

Under these conditions, the series $U(a,b)$ is the same as in case 2 of the section I (see (I.8). For even n, $U(a,b)$ is not divisible by n, but for odd n, it is proportional to $n^2$.

On the other hand, because (a+b) and c are both divisible by n, the series $U((a+b),c)$ is proportional to $n^n$. Therefore, the divisibility of $U(a,b,c)$ is determined by the divisibility of $U(a,b)$. Thus we have

$$[U(a, b, c), n^2] = 1 \quad \text{if n is odd} \quad (II.4a)$$

$$[U(a, b, c), n] = 0 \quad \text{if n is even} \quad (II.4b)$$

## Conclusion

Solutions to problems often require the clarification of the compatibility of different terms of an equation with respect to evenness, co-prime status or divisibility. The particular truncated binomial series described in this paper are related to equations of type $An^k - U(a, b, c, d ...) = 0$, where $U(a,b,c,d ....)$ is a TBS. The choice of $k$ depends on the divisibility property of the series $U(a,b,c,d)$.